\documentclass[11pt]{article} \makeatother
\usepackage{amsfonts, amsmath, amssymb}
\usepackage{graphicx}
\usepackage{subfigure}
\allowdisplaybreaks[4]
\newtheorem{theo}{Theorem}[section]
\newtheorem{lemma}[theo]{Lemma}
\newtheorem{proposition}[theo]{Proposition}

\newcommand{\be}{\begin{equation}}
\newcommand{\ee}{\end{equation}}
\newcommand\bes{\begin{eqnarray}} \newcommand\ees{\end{eqnarray}}
\newcommand{\bess}{\begin{eqnarray*}}
\newcommand{\eess}{\end{eqnarray*}}

\numberwithin{equation}{section}

\begin{document}
\date{}
\setlength{\baselineskip}{17pt}{\setlength\arraycolsep{2pt}

\font\sevenrm=cmr7
\font\sevenit=cmti7
\font\sevenbf=cmbx7
\makeatletter
\def\ps@plain{
\def\@oddhead{\ifnum\thepage=1\hss\baselineskip8pt
\vtop to 0 pt{\vskip-1truecm\hbox%{\sevenit
%Electronic Journal of Qualitative Theory of Differential Equations}\hbox{\sevenrm 2010, No. {\sevenbf 15}, 1-12;
%{\sevenbf http://www.math.u-szeged.hu/ejqtde/}
}\vss}\else\hss\fi}\def\@oddfoot{\rm\hfil EJQTDE, 2010 No. 15, p. \thepage}
\def\@evenhead{}\let\@evenfoot\@oddfoot}
\makeatother
\pagestyle{plain}
\thispagestyle{plain}

\begin{center} {\bf\Large General decay  for a viscoelastic wave equation with dynamic boundary conditions and a time-varying delay }\\[2mm]

 {\large Gang Li, Biqing Zhu, Danhua Wang  \\[1mm]
{College of Mathematics and Statistics, Nanjing University of
Information Science and Technology, Nanjing 210044, China.\\ E-mail: {\bf brucechu@163.com}.}}\\[1mm]
%{Department of Mathematics,  Southeast University, Nanjing
%210096, China.}\\[12mm]
\end{center}

\begin{abstract}
The goal of this paper is to study a nonlinear viscoelastic wave equation with strong damping, time-varying delay and dynamical boundary condition. By introducing suitable energy and Lyapunov functionals, under suitable assumptions, we then prove a general decay result of the energy, from which the usual exponential and polynomial decay rates are only special cases.
\end{abstract}



\section{Introduction }
\setcounter{equation}{0}

In this paper, we consider the following problem:
\begin{equation}\label{1.1}
\left\{ {{\begin{array}{*{20}l} \displaystyle u_{tt}-\Delta u +\int_{0}^{t}g(t-s)\Delta u(x,s){\rm d}s -\alpha\Delta u_{t} =0 , \ & x\in \Omega,  t > 0, \medskip\\
\displaystyle u(x,t) =0, \ & x\in \Gamma_{0},  t > 0, \medskip\\
\displaystyle u_{tt}(x, t) = \displaystyle-\frac {\partial u}{\partial \nu}(x,t)+\int_{0}^{t}g(t-s)\frac{\partial u}{\partial \nu}(x,s){\rm d}s-\alpha\frac {\partial u_{t}}{\partial \nu}(x,t)\medskip\\\quad\quad\quad\quad\ \ -\mu_1u_t(x,t)-\mu_2u_t(x,t-\tau(t)), \ & x \in \Gamma_{1}, t > 0, \medskip\\
u(x, 0) = u_0(x), u_t(x, 0) = u_1(x), \ & x \in  \Omega,\medskip\\
u_t(x,t-\tau(0))=f_0(x,t-\tau(0)),\ & x \in \Gamma_{1}, t \in(0,\tau(0)),
\end{array} }} \right.
\end{equation}
where  $\Omega$ is a regular and bounded domain of $\mathbb{R}^{N}$, $(N\geq1) $, $\partial \Omega= \Gamma_{0}\cup\Gamma_1$, $mes(\Gamma_0)>0$, ${\Gamma_0\cap\Gamma_1= \emptyset}$ and $\frac {\partial}{\partial \nu}$ denotes the unit outer normal derivative. Moreover, $\tau(t)>0$ is the time-varying delay term, $\alpha$, $\mu_1$ and $\mu_2$ are positive constants. The initial datum $u_0$, $u_1$  and $f_0$ are given functions belonging to suitable spaces.

From the mathematical point of view, these problems like \eqref{1.1} take into account acceleration terms on the boundary. Such type of boundary conditions are usually called dynamic boundary conditions (see  \cite{AKS1996}, \cite{BEA1976}, \cite{BST1964}, \cite{GOL2006} for more details). The above model without delay term (i.e., $\mu_2=0$), has been studied by many authors in recent years. For example, Gerbi and Said-Houair in \cite{GER2008} studied problem \eqref{1.1} with source term $|u|^{p-2}u$ and nonlinear damping on the boundary but without the relaxation function $g$. They showed  that if the initial data are large enough then the energy and the $L^p$ norm of the solution of the problem is unbounded, grows up exponentially as time goes to infinity. Later in \cite{GH2011}, they established the global existence and asymptotic stability of solutions starting in a stable set by combining the potential well method and the energy method. A blow-up result for the case $m=2$ with initial data  in the unstable set was also obtained. Recently, when the relaxation function $g\neq 0$, they in \cite{GH2013} got the existence and exponential growth results. For the other works, we refer the readers to (\cite{GOL2006}, \cite{GH2012}, \cite{SW2005}, \cite{RU2010})  and the references therein.

On the other hand, the above model with delay term (i.e., $\mu_2\neq0$) has become an active area of research. The delay term may be a source of instability, we refer the readers to  (\cite{BEN2014}, \cite{NIC2006}, \cite{NIC2008}, \cite{WU2013})  and the references therein. For example, in \cite{GER2012}, Stephane Gerbi and Belkacem Said-Houari considered the following linear damped wave equation with dynamic boundary conditions and a delay boundary term:
\begin{equation}
\left\{ {{\begin{array}{*{20}l} \displaystyle u_{tt}-\Delta u -\alpha\Delta u_{t} =0 , \ & x\in \Omega,  t > 0, \medskip\\
\displaystyle u(x,t) =0, \ & x\in \Gamma_{0},  t > 0, \medskip\\
\displaystyle u_{tt}(x, t) = \displaystyle-\left(\frac {\partial u}{\partial \nu}(x,t)+\alpha\frac {\partial u_{t}}{\partial \nu}(x,t) +\mu_1u_t(x,t)+\mu_2u_t(x,t-\tau)\right), \ & x \in \Gamma_{1}, t > 0, \medskip\\
u(x, 0) = u_0(x), u_t(x, 0) = u_1(x), \ & x \in  \Omega,\medskip\\
u_t(x,t-\tau)=f_0(x,t-\tau))\ & x \in \Gamma_{1}, t \in(0,\tau),
\end{array} }} \right.\nonumber
\end{equation}
 under the condition that if the weight of the delay term in the feedback is less than the weight of the term without delay or if it is greater under an assumption between the damping factor and the difference of two weights, they proved the global existence of the solutions and the exponential stability of the system.
Later, Mohamed FERHAT and Ali HAKEM in \cite{FER2015} considered the following wave equation with dynamic boundary conditions:
\begin{equation}
\left\{ {{\begin{array}{*{20}l} \displaystyle u_{tt}-\Delta u -\alpha\Delta u_{t} -\int_{0}^{t}g(t-s)\Delta u(x,s){\rm d}s=|u|^{p-1}u ,\ & {\rm in}\ \  \Omega\times(0,+\infty), \medskip\\
\displaystyle u(x,t) =0, \ & {\rm on}\ \ \Gamma_{0}\times(0,+\infty),  \medskip\\
\displaystyle u_{tt}(x, t) = \displaystyle-a\left[\frac {\partial u}{\partial \nu}(x,t)-\alpha\frac {\partial u_{t}}{\partial \nu}(x,t)-\int_{0}^{t}g(t-s)\frac{\partial u}{\partial \nu}(x,s){\rm d}s\right.\medskip\\\quad\quad\quad\quad\ \ \bigg.+\mu_1\psi(u_t(x,t))+\mu_2\psi(u_t(x,t-\tau))\bigg], \ & {\rm on}\ \ \Gamma_{1}\times(0,+\infty), \medskip\\
u(x, 0) = u_0(x), u_t(x, 0) = u_1(x), \ & x \in  \Omega,\medskip\\
u_t(x,t-\tau)=f_0(x,t-\tau),\ & {\rm on}\ \ \Gamma_{1}\times(0,+\infty).\nonumber
\end{array} }} \right.
\end{equation}
By using the potential well method and introducing suitable Lyapunov function, they proved the global existence and established general decay estimates for the energy.

Recently, the case of time-varying delay has been studied by (\cite{BEN2013}, \cite{NIC2011}, \cite{NIC2009}).  For example, Nicaise, Valein and Fridman \cite{NIC2009} in one space dimension. They proved the exponential stability result under the condition
\begin{equation}\label{1.2}
\mu_2<\sqrt{1-d}\mu_1
\end{equation}
where $d$ is a constant such that
\begin{equation}\label{1.3}
\tau'(t)\leq d<1, \quad \forall t>0.
\end{equation}
Later, Serge Nicaise, Cristina Pignotti and Julie Valein considered the following problem
\begin{equation}\label{1.4}
\left\{ {{\begin{array}{*{20}l} \displaystyle u_{tt}-\Delta u  =0 , \ & {\rm in}\ \  \Omega\times(0,\infty), \medskip\\
\displaystyle u(x,t) =0, \ & {\rm on}\ \ \Gamma_{D}\times(0,\infty), \medskip\\
\displaystyle \displaystyle\frac {\partial u}{\partial \nu}(x,t)=-\mu_1u_t(x,t)-\mu_2u_t(x,t-\tau(t)), \ & {\rm on}\ \ \Gamma_{N}\times(0,\infty), \medskip\\
u(x, 0) = u_0(x), u_t(x, 0) = u_1(x), \ & {\rm in}\ \  \Omega,\medskip\\
u_t(x,t-\tau(0))=f_0(x,t-\tau(0)),\ & {\rm on}\ \ \Gamma_{N}\times(0,\tau(0)),
\end{array} }} \right.
\end{equation}
they extend the last result to general space dimension under the hypothesis
\begin{equation}
\tau(t)\geq \tau_0>0, \quad \forall t>0,\nonumber
\end{equation}
assumed in \cite{NIC2009}, that is the delay may degenerate. They also gave a well-posedness result and an exponential stability estimate for problem \eqref{1.4} under
a suitable relation between the coefficients.

Motivated by these results, in this paper, we intend to study the general decay result to problem \eqref{1.1}.
%The main difficulties we encounter here arise from the simultaneous appearance of the viscoelastic term, the time-varying delay term and the nonlinear boundary damping term.
 Our main contribution is an extension of previous result from \cite{GER2012} to relaxation function $g$ and time-varying delays with $\tau(t)\geq 0$. By introducing new energy and Lyapunov functionals, we show in this article that the decay rates of the solution energy is similar to the relaxation function, which are not necessarily decaying like polynomial or exponential functions.
%We use the Faedo-Galerkin approximation, potential well theory, perturbed energy method and concavity technique to get the above results.

The paper is organized as follows. In Section 2, we present some assumptions needed for our work and state the main result. The general decay result is given in Section 3.

\section{Preliminaries }\label{2}
\setcounter{equation}{0}
In this section we present some assumptions and state the main result. For the relaxation function $g$, we assume the following

$(G1)$ $g$: $\mathbb{R}_{+}\longrightarrow \mathbb{R}_{+}$ is a nonincreasing differentiable function satisfying
$$g(0)>0,\quad 1-\int_{0}^{\infty}g(s){\rm d}s=l>0.$$

$(G2)$ There exists a nonincreasing differentiable function $\xi: \mathbb{R}_{+}\longrightarrow \mathbb{R}_{+}$ such that
$$g'(s)\leq-\xi(s) g(s),\quad \forall s \in\mathbb{R}_{+}$$
and
$$\int_0^{+\infty}\xi(t){\rm d}t=\infty.$$

We denote $H_{\Gamma_0}^{1}=\{u\in H^{1}(\Omega)|u_{|_{\Gamma_0}}=0\}$, $\mathcal{V}=H_{\Gamma_{0}}^{1}(\Omega)\cap L^{2}(\Gamma_{1})$ and by $(\cdot,\cdot)$ we denote the scalar product in $L^{2}(\Omega)$; i.e.,
$$(u,v)(t)=\int_{\Omega}u(x,t)v(x,t){\rm d}x.$$
%Also, by $\|\cdot\|_{q}$ we mean the $L^{q}(\Omega)$ norm for $1\leq q\leq\infty$, by $\|\cdot\|_{q,\Gamma_{1}}$ the $L^{q}(\Gamma_{1})$ norm.
%We will also use the embedding (see \cite{ADA1975}):
%$H_{\Gamma_0}^{1}\hookrightarrow L^{q}(\Gamma_1)$,\quad $2\leq q\leq \overline q$, where $\overline q$ satifies
%\begin{align}
%\overline q =\left\{ {{\begin{array}{*{20}l} \frac {2(N-1)}{N-2} > 1, \ & if \ \ N \geq3, \medskip\\ +\infty, \ & if \ \ N=1,2. \end{array} }} \right.\nonumber
%\end{align}
As in \cite{NIC2006}, let us introduce the new variable
\begin{align}
z(x,\rho,t)=u_t(x,t-\tau(t)\rho), \quad x\in\Gamma_1, \rho\in(0,1), t>0.\nonumber
\end{align}
Then, we have
\begin{align}
\tau(t)z_t(x,\rho,t)+z_\rho(x,\rho,t)=0, \quad {\rm in} \ \ \Gamma_1\times(0,1)\times(0,+\infty).\nonumber
\end{align}
Therefore, problem \eqref{1.1} is equivalent to
\begin{equation}\label{2.1}
\left\{ {{\begin{array}{*{20}l} \displaystyle u_{tt}-\Delta u +\int_{0}^{t}g(t-s)\Delta u(x,s){\rm d}s -\alpha\Delta u_{t} =0 , \ & x\in \Omega,  t > 0, \medskip\\
\displaystyle \tau(t)z_t(x,\rho,t)+z_\rho(x,\rho,t)=0,  \ &  x\in\Gamma_1, \rho\in(0,1), t>0, \medskip\\
\displaystyle u(x,t) =0, \ & x\in \Gamma_{0},  t > 0, \medskip\\
\displaystyle u_{tt}(x, t) = \displaystyle-\frac {\partial u}{\partial \nu}(x,t)+\int_{0}^{t}g(t-s)\frac{\partial u}{\partial \nu}(x,s){\rm d}s\medskip\\ \quad\quad\quad\quad\ \ \displaystyle -\alpha\frac {\partial u_{t}}{\partial \nu}(x,t) -\mu_1u_t(x,t)-\mu_2z(x,1,t), \ & x \in \Gamma_{1}, t > 0, \medskip\\
z(x,0,t)=u_t(x,t),\ & x \in \Gamma_{1}, t > 0, \medskip\\
u(x, 0) = u_0(x), u_t(x, 0) = u_1(x), \ & x \in  \Omega,\medskip\\
z(x,\rho,0)=f_0(x,-\rho\tau(0)),\ & x \in \Gamma_{1}, \rho \in(0,1).
\end{array} }} \right.
\end{equation}

We now state, without a proof, local existence result, which can be established by using the Fadeo-Galerkin approximation method (see \cite{NIC2011}, \cite{NIC2009} for more details).

\begin{lemma}\label{Lemma 2.1}
  Suppose that $(G1)$ and $(G2)$ hold. Then given $u_0\in H_{\Gamma_0}^{1}(\Omega)$, $u_{1}\in L^{2}(\Omega)$ and $f_0\in L^2(\Omega\times(0,1))$, then there exist $T>0$ and a unique weak solution $(u,z)$ of problem \eqref{2.1} on $(0,T)$ satisfying
$$u\in C\left([0,T],H_{\Gamma_0}^{1}(\Omega)\right)\cap C^{1}\left([0,T],L^{2}(\Omega)\right),$$
$$u_{t}\in L^2\left(0,T;H_{\Gamma_0}^{1}(\Omega)\right)\cap L^{2}\left((0,T)\times \Gamma_1\right).$$
\end{lemma}

We define the new energy of system \eqref{1.1} as
\begin{align}\label{2.2}
E(t):=&\frac{1}{2}\left[\|u_t(t)\|_{2}^{2}+\left(1-\int_{0}^{t}g(s){\rm d}s\right)\|\nabla u(t)\|_{2}^{2}+(g\circ \nabla u)(t)+\|u_t(t)\|_{2,\Gamma_1}^{2}\right]\nonumber\\
&+\frac{\zeta}{2}\tau(t)\int_{0}^{1}\int_{\Gamma_1}u_{t}^{2}(x,t-\tau(t)\rho){\rm d}\rho{\rm d}\sigma,
\end {align}
where $\zeta$ is a positive constant such that
\begin{align}\label{2.3}
2\mu_1-\frac{\mu_2}{\sqrt{1-d}}-\zeta>0 \quad {\rm and} \quad \zeta-\frac{\mu_2}{\sqrt{1-d}} >0
\end {align}
and
$$(g\circ\nabla u)(t)=\int_{0}^{t}g(t-s)\|\nabla u(t)-\nabla u(s)\|_{2}^{2}ds\geq0.$$
Then, we state the main result as follows
\begin{theo}\label{th5.3}
Let $(u_0,u_1)\in H_{\Gamma_0}^1\times L^{2}(\Omega)$ be given. Assume that $g$ and $\xi$ satisfy $(G1)$ and $(G2)$. Then, for each $t_0>0$, there exist two positive constants $K$ and $k$ such that, for any solution of the problem \eqref{1.1}, the energy satisfies
\begin{align}\label{5.13}
E(t)\leq Ke^{-k\int_{t_0}^t\xi(s){\rm d}s}.
\end {align}
\end{theo}

\section{Decay of solutions} \label{5}

As mentioned earlier, in this section, we prove the general decay result for problem \eqref{1.1} under the assumption \eqref{1.2}.
\begin{proposition}\label{Proposition 5.1}
For any regular solution of problem \eqref{1.1} we have
\begin{align}\label{5.3}
E'(t)=&-\alpha\|\nabla u_t(t)\|_{2}^{2}-\mu_1\|u_t(t)\|_{2,\Gamma_1}^2-\mu_2\int_{\Gamma_1}u_t(x,t)u_{t}(x,t-\tau(t)){\rm d}\sigma+\frac{1}{2}(g'\circ\nabla u)(t)\nonumber\\
&-\frac{1}{2}g(t)\|\nabla u(t)\|_{2}^{2}-\frac{\zeta}{2}\int_{\Gamma_1}u_{t}^{2}(x,t-\tau(t))(1-\tau'(t)){\rm d}\sigma+\frac{\zeta}{2}\|u_t(t)\|_{2,\Gamma_1}^{2}.
\end {align}
\end{proposition}
{\bf Proof.} Differentiating \eqref{2.2} we get
\begin{align}\label{5.4}
E'(t)=&\int_{\Omega}u_{tt}u_{t}{\rm d}x+\left(1-\int_{0}^{t}g(s){\rm d}s\right)\int_{\Omega}\nabla u\nabla u_t{\rm d}x-\frac{1}{2}g(t)\|\nabla u\|_{2}^{2}+\frac{1}{2}\left(g'\circ\nabla u\right)(t)\nonumber\\
&+\int_{\Gamma_1}u_{tt}u_{t}{\rm d}\sigma+\frac{\zeta}{2}\tau'(t)\int_{0}^{t}\int_{\Gamma_1}u_{t}^{2}(x,t-\tau(t)\rho){\rm d}\rho{\rm d}\sigma\nonumber\\
&+\zeta\tau(t)\int_{0}^{t}\int_{\Gamma_1}u_{tt}(x,t-\tau(t)\rho)u_{t}(x,t-\tau(t)\rho)(1-\tau'(t)\rho){\rm d}\rho{\rm d}\sigma.
\end {align}

{\bf Case 1.} If $\tau(t)\neq 0$, then
\begin{align}
u_{t}(x,t-\tau(t)\rho)=-\tau^{-1}(t)u_\rho(x,t-\tau(t)\rho)\nonumber
\end {align}
and
\begin{align}
u_{tt}(x,t-\tau(t)\rho)=\tau^{-2}(t)u_{\rho\rho}(x,t-\tau(t)\rho).\nonumber
\end {align}
So we get
\begin{align}\label{5.5}
&\int_{0}^{1}u_{tt}(x,t-\tau(t)\rho)u_{t}(x,t-\tau(t)\rho)(1-\tau'(t)\rho){\rm d}\rho\nonumber\\
=&-\tau^{-3}(t)\int_{0}^{1}u_{\rho\rho}(x,t-\tau(t)\rho)u_{\rho}(x,t-\tau(t)\rho)(1-\tau'(t)\rho){\rm d}\rho\nonumber\\
=&-\tau^{-3}(t)\left[u_{\rho}^{2}(x,t-\tau(t)\rho)(1-\tau'(t)\rho)\right]_{0}^{1}\nonumber\\
&+\tau^{-3}(t)\int_{0}^{1}u_{\rho\rho}(x,t-\tau(t)\rho)u_{\rho}(x,t-\tau(t)\rho)(1-\tau'(t)\rho){\rm d}\rho\nonumber\\
&-\tau'(t)\tau^{-3}(t)\int_{0}^{1}u_{\rho}(x,t-\tau(t)\rho)u_{\rho}(x,t-\tau(t)\rho){\rm d}\rho\nonumber\\
=&-\frac{1}{2}\tau'(t)\tau^{-3}(t)\int_{0}^{t}u_{\rho}^{2}(t-\tau(t)\rho){\rm d}\rho\nonumber\\
&-\frac{\tau^{-1}(t)}{2}u_{t}^{2}(x,t-\tau(t))(1-\tau'(t))+\frac{\tau^{-1}(t)}{2}u_{t}^{2}(x,t)\nonumber\\
=&-\frac{1}{2}\tau'(t)\tau^{-1}(t)\int_{0}^{t}u_{t}^{2}(x,t-\tau(t)\rho){\rm d}\rho\nonumber\\
&-\frac{\tau^{-1}(t)}{2}u_{t}^{2}(x,t-\tau(t))(1-\tau'(t))+\frac{\tau^{-1}(t)}{2}u_{t}^{2}(x,t).
\end {align}
By using \eqref{5.4}, \eqref{5.5} and the boundary condition on $\Gamma_1$, we obtain \eqref{5.3}.

{\bf Case 2.} If $\tau(t)= 0$, then from \eqref{5.4}, we get
\begin{align}\label{5.6}
E'(t)=&-\alpha\|\nabla u_t(t)\|_{2}^{2}-(\mu_1+\mu_2)\|u_t(t)\|_{2,\Gamma_1}^2+\frac{1}{2}(g'\circ\nabla u)(t)-\frac{1}{2}g(t)\|\nabla u(t)\|_{2}^{2}\nonumber\\
&+\frac{\zeta}{2}\|u_t(t)\|_{2,\Gamma_1}^{2}.
\end {align}
Therefore, \eqref{5.3} is proved for all times $t>0$.
\begin{lemma}\label{Lemma 5.2}
For any regular solution of problem \eqref{1.1} the energy decays and there exists a positive constant $C$ such that
\begin{align}\label{5.7}
E'(t)\leq&-\alpha\|\nabla u_t(t)\|_{2}^{2}+\frac{1}{2}(g'\circ\nabla u)(t)-\frac{1}{2}g(t)\|\nabla u(t)\|_{2}^{2}\nonumber\\
&-C\left(\|u_t(t)\|_{2,\Gamma_1}^{2}+\|u_t(t-\tau(t))\|_{2,\Gamma_1}^{2}\right).
\end {align}
\end{lemma}
{\bf Proof.} In the case of $\tau(t)\neq0$, by Cauchy-Schwarz's inequality, we have
\begin{align}
E'(t)\leq&-\alpha\|\nabla u_t(t)\|_{2}^{2}+\frac{1}{2}(g'\circ\nabla u)(t)-\frac{1}{2}g(t)\|\nabla u(t)\|_{2}^{2}-\mu_1\|u_t(t)\|_{2,\Gamma_1}^{2}+\frac{\mu_2}{2\sqrt{1-d}}\|u_t(t)\|_{2,\Gamma_1}^{2}\nonumber\\
&+\frac{\mu_2\sqrt{1-d}}{2}\|u_t(x,t-\tau(t))\|_{2,\Gamma_1}^{2}-\frac{\zeta}{2}(1-\tau'(t))\|u_t(x,t-\tau(t))\|_{2,\Gamma_1}^{2}+\frac{\zeta}{2}\|u_t(t)\|_{2,\Gamma_1}^{2}\nonumber\\
\leq&-\alpha\|\nabla u_t(t)\|_{2}^{2}+\frac{1}{2}(g'\circ\nabla u)(t)-\frac{1}{2}g(t)\|\nabla u(t)\|_{2}^{2}-C\left(\|u_t(t)\|_{2,\Gamma_1}^{2}+\|u_t(t-\tau(t))\|_{2,\Gamma_1}^{2}\right)\nonumber,
\end {align}
by \eqref{2.3} we easily get  \eqref{5.7}. In the case of $\tau(t)=0$, when $\displaystyle \zeta<2\mu_1<\frac{2(\mu_1+\mu_2)}{d}$, by \eqref{5.6} we obtain \eqref{5.7}.

Now, we use the following modified functional, for positive constants $\varepsilon_1$, $\varepsilon_2$ and $\varepsilon_3$, we have
\begin{align}\label{5.8}
L(t)=E(t)+\varepsilon_1 \psi(t)+\varepsilon_2 \phi(t)+\varepsilon_3 I(t),
\end {align}
where
\begin{align}\label{5.9}
\psi(t)=\int_\Omega u_tu{\rm d}x+\int_{\Gamma_1} u_tu{\rm d}\sigma+\frac{\alpha}{2}\|\nabla u\|_2^2,
\end {align}
\begin{align}\label{5.10}
\phi(t)=-\int_\Omega u_t\int_0^{t}g(t-s)\left(u(t)-u(s)\right){\rm d}s{\rm d}x
\end {align}
and
\begin{align}\label{5.11}
I(t)=-\zeta\tau(t)\int_{\Gamma_1}\int_{0}^{1}e^{-2\tau(t)\rho}u_{t}^{2}(x,t-\tau(t)\rho){\rm d}\rho{\rm d}\sigma.
\end {align}
It is easy to check that, by using Poincare's inequality, trace inequality, \eqref{2.3} and for $\varepsilon_1$, $\varepsilon_2$, $\varepsilon_3$ small enough, there exists two constants $\alpha_1$ and $\alpha_2$ such that
\begin{align}\label{5.12}
\alpha_1L(t)\leq E(t)\leq \alpha_2L(t).
\end {align}

Next, we estimate the derivative of $L(t)$ according to the following lemmas.
\begin{lemma}\label{Lemma 5.4}
Under the conditions of Theorem \ref{th5.3}, the functional $\psi(t)$ defined in \eqref{5.9} satisfies
\begin{align}\label{5.14}
\psi'(t)\leq &\|u_{t}(t)\|_{2}^{2}+\left(1+\frac{\mu_{1}^{2}}{4\delta}\right)\| u_t(t)\|_{2,\Gamma_1}^2+\left(\delta-l+2c\delta \right)\|\nabla u(t)\|_2^2\nonumber\\
&+\frac{\mu_{2}^2}{4\delta}\|u_t(t-\tau(t))\|_{2,\Gamma_1}^{2}+\frac{1-l}{4\delta}\left(g\circ\nabla u\right)(t),
\end {align}
for some $\delta>0$.
\end{lemma}
{\bf Proof.} By using the differential equation in \eqref{1.1}, we get
\begin{align}\label{5.15}
\psi'(t)=&\| u_t(t)\|_{2}^{2}+\int_\Omega u_{tt}(t)u(t){\rm d}x+\int_{\Gamma_1} u_{tt}(t)u(t){\rm d}\sigma+\|u_{t}(t)\|_{2,\Gamma_1}^{2}+\alpha\int_\Omega \nabla u_t(t)\cdot\nabla u(t){\rm d}x\nonumber\\
=&\| u_t(t)\|_{2}^{2}+\|u_{t}(t)\|_{2,\Gamma_1}^{2}+\int_{\Gamma_1}\frac{\partial u(t)}{\partial \nu}u(t){\rm d}\sigma-\| \nabla u(t)\|_{2}^{2}-\int_{0}^{t}g(t-s)\int_{\Gamma_1}\frac{\partial u(s)}{\partial \nu}u(s){\rm d}\sigma{\rm d}s\nonumber\\
&+\int_{0}^{t}g(t-s)\int_\Omega\nabla u(s)\cdot\nabla u(t){\rm d}x{\rm d}s+\alpha\int_{\Gamma_1}\frac{\partial u_t(t)}{\partial \nu}u(t){\rm d}\sigma-\mu_1\int_{\Gamma_1}u_t(t)u(t){\rm d}\sigma\nonumber\\
&-\mu_2\int_{\Gamma_1}u_t(t-\tau(t))u(t){\rm d}\sigma\nonumber\\
=&\| u_t(t)\|_2^2+\|u_{t}(t)\|_{2,\Gamma_1}^{2}-\|\nabla u(t)\|_2^2+\int_\Omega\nabla u(t)\cdot\int_0^tg(t-s)\nabla u(s){\rm d}s{\rm d}x\nonumber\\
&-\mu_1\int_{\Gamma_1}u_t(t)u(t){\rm d}\sigma-\mu_2\int_{\Gamma_1}u_t(t-\tau(t))u(t){\rm d}\sigma.
\end {align}
We now estimate the right hand side of \eqref{5.15}. For a positive constant $\delta$, we have the estimates as follows
\begin{align}\label{5.16}
\int_\Omega\nabla u(t)\cdot\int_0^tg(t-s)\nabla u(s){\rm d}s{\rm d}x\leq(\delta+1-l)\|\nabla u(t)\|_2^2+ \frac{1-l}{4\delta}(g\circ \nabla u)(t).
\end {align}
By Young's inequality and trace inequality, we have
\begin{align}\label{5.17}
-\mu_1\int_{\Gamma_1}u_t(t)u(t){\rm d}\sigma-\mu_2\int_{\Gamma_1}u_t(t-\tau(t))u(t){\rm d}\sigma\leq& 2c\delta\|\nabla u(t)\|_{2}^{2}+\frac{\mu_{1}^{2}}{4\delta}\|u_{t}\|_{2,\Gamma_1}^{2}\nonumber\\
&+\frac{\mu_{2}^{2}}{4\delta}\|u_t{(t-\tau(t))}\|_{2,\Gamma_1}^{2}.
\end {align}
Combining \eqref{5.15}-\eqref{5.17}, we arrive at \eqref{5.14}.
\begin{lemma}\label{Lemma 5.5}
Under the conditions of Theorem \ref{th5.3}, the functional $\phi(t)$ defined in \eqref{5.10} satisfies
\begin{align}\label{5.18}
\phi'(t)\leq&\left(\delta-g(0)\right)\|u_t(t)\|_2^2+\mu_1\|u_{t}(t)\|_{2,\Gamma_1}^{2}+\left[\delta+2\delta(1-l)^2\right]\|\nabla u(t)\|_2^2+\frac{\alpha}{2}\|\nabla u_{t}(t)\|_{2}^{2}\nonumber\\
&+\left[\frac{1-l}{4\delta}+\left(2\delta+\frac{1}{4\delta}\right)(1-l)^2+\frac{\alpha(1-l)}{2}+\frac{1-l}{2\delta\lambda_1}\right](g\circ\nabla u)(t)+\frac{1-l}{4\delta\lambda_1}(-g'\circ\nabla u)(t)\nonumber\\
&+\mu_2\|u_t{(t-\tau(t))}\|_{2,\Gamma_1}^{2},
\end {align}
for some $\delta>0$.
\end{lemma}
{\bf Proof.} By using the differential equation in \eqref{1.1}, we get
\begin{align}\label{5.19}
\phi'(t)=&-\int_\Omega u_{tt}(t)\int_0^{t}g(t-s)\left(u(t)-u(s)\right){\rm d}s{\rm d}x-\int_\Omega u_{t}(t)\int_0^{t}g'(t-s)\left(u(t)-u(s)\right){\rm d}s{\rm d}x\nonumber\\
&-\left(\int_0^{t}g(s){\rm d}s\right)\int_\Omega|u_t(t)|^{2}{\rm d}x\nonumber\\
=&\int_\Omega\nabla u(t)\cdot\int_0^{t}g(t-s)\left(\nabla u(t)-\nabla u(s)\right){\rm d}s{\rm d}x-\left(\int_0^{t}g(s){\rm d}s\right)\int_\Omega u_{t}^2(t){\rm d}x\nonumber\\
&-\int_\Omega\left(\int_{0}^{t}g(t-s)\nabla u(s){\rm d} s\right)\left(\int_0^{t}g(t-s)\left(\nabla u(t)-\nabla u(s)\right){\rm d}s\right){\rm d}x\nonumber\\
&+\alpha\int_\Omega\nabla u_t(t)\cdot\int_0^{t}g(t-s)\left(\nabla u(t)-\nabla u(s)\right){\rm d}s{\rm d}x\nonumber\\
&-\int_\Omega u_t(t)\int_0^{t}g'(t-s)\left(u(t)-u(s)\right){\rm d}s{\rm d}x\nonumber\\
&+\mu_1\int_{\Gamma_1}u_t(t)\int_0^{t}g(t-s)\left(u(t)-u(s)\right){\rm d}s{\rm d}\sigma\nonumber\\
&+\mu_2\int_{\Gamma_1}u_t(t-\tau(t))\int_0^{t}g(t-s)\left(u(t)-u(s)\right){\rm d}s{\rm d}\sigma.
\end {align}
We now estimate the right side of \eqref{5.19}, using Young's inequality, Hoider's inequality and  Cauchy-Schwarzs's inequality, we get
\begin{align}\label{5.20}
\int_\Omega\nabla u(t)\cdot\int_0^{t}g(t-s)\left(\nabla u(t)-\nabla u(s)\right){\rm d}s{\rm d}x\leq&\delta\|\nabla u(t)\|_2^2+\frac{1-l}{4\delta}(g\circ\nabla u)(t),
\end {align}
\begin{align}\label{5.21}
&-\int_\Omega\left(\int_{0}^{t}g(t-s)\nabla u(s){\rm d} s\right)\left(\int_0^{t}g(t-s)\left(\nabla u(t)-\nabla u(s)\right){\rm d}s\right){\rm d}x\nonumber\\
\leq&\delta\int_\Omega\left(\int_{0}^{t}g(t-s)\nabla u(s){\rm d} s\right)^{2}{\rm d}x+\frac{1}{\delta}\int_\Omega\left(\int_0^{t}g(t-s)\left(\nabla u(t)-\nabla u(s)\right){\rm d}s\right)^{2}{\rm d}x\nonumber\\
\leq&2\delta(1-l)^{2}\|\nabla u(t)\|_2^{2}+\left(2\delta+\frac{1}{4\delta}\right)(1-l)^{2}(g\circ\nabla u)(t),
\end {align}
\begin{align}\label{5.22}
&\alpha\int_\Omega\nabla u_t(t)\cdot\int_0^{t}g(t-s)\left(\nabla u(t)-\nabla u(s)\right){\rm d}s{\rm d}x\nonumber\\
\leq&\frac{\alpha}{2}\|\nabla u_t(t)\|_2^2+\frac{\alpha(1-l)}{2}(g\circ\nabla u)(t),
\end {align}
\begin{align}\label{5.23}
&\mu_1\int_{\Gamma_1}u_t(t)\int_0^{t}g(t-s)\left(u(t)-u(s)\right){\rm d}s{\rm d}\sigma\nonumber\\
\leq&\mu_1\|u_t(t)\|_{2,\Gamma_1}^{2}+\frac{1-l}{4\delta\lambda_1}(g\circ\nabla u)(t),
\end {align}
\begin{align}\label{5.24}
&\mu_2\int_{\Gamma_1}u_t(t-\tau(t))\int_0^{t}g(t-s)\left(u(t)-u(s)\right){\rm d}s{\rm d}\sigma\nonumber\\
\leq&\mu_2\|u_t(t-\tau(t))\|_{2,\Gamma_1}^{2}+\frac{1-l}{4\delta\lambda_1}(g\circ\nabla u)(t),
\end {align}
since $g$ is continuous and $g(0)>0$, then for any $t_0>0$, we have
\begin{align}\label{5.25}
\int_{0}^{t}g(s){\rm d} s\geq\int_{0}^{t_0}g(s){\rm d} s=g_0, \quad \forall t\geq t_0,
\end {align}
then we use \eqref{5.25} to get
\begin{align}\label{5.26}
&-\int_\Omega u_t(t)\int_0^{t}g'(t-s)\left(u(t)-u(s)\right){\rm d}s{\rm d}x-\left(\int_0^{t}g(s){\rm d}s\right)\int_\Omega u_{t}^2(t){\rm d}x\nonumber\\
\leq&\delta\| u_t(t)\|_{2}^{2}+\frac{1-l}{4\delta\lambda_1}\left(-g'\circ\nabla u\right)(t)-g_0\| u_t(t)\|_{2}^{2}.
\end {align}
A combination of \eqref{5.19}-\eqref{5.26} yields \eqref{5.18}.
\begin{lemma}\label{Lemma 5.6}
Under the conditions of Theorem \ref{th5.3}, the functional $I(t)$ defined in \eqref{5.11} satisfies
\begin{align}\label{5.27}
I'(t)\leq-2I(t)+\zeta\|u_t(t)\|_{2,\Gamma_1}^{2}
\end {align}
for some $\delta>0$.
\end{lemma}
{\bf Proof.} Differentiating \eqref{5.11} we have
\begin{align}\label{5.28}
I'(t)=&\zeta\tau'(t)\int_{\Gamma_1}\int_{0}^{1}e^{-2\tau(t)\rho}u_{t}^{2}(x,t-\tau(t)\rho){\rm d}\rho{\rm d}\sigma\nonumber\\
&-2\zeta\tau'(t)\tau(t)\int_{\Gamma_1}\int_{0}^{1}e^{-2\tau(t)\rho}u_{t}^{2}(x,t-\tau(t)\rho){\rm d}\rho{\rm d}\sigma\nonumber\\
&+2\zeta\tau(t)\int_{\Gamma_1}\int_{0}^{1}e^{-2\tau(t)\rho}u_{t}(x,t-\tau(t)\rho)u_{tt}(x,t-\tau(t)\rho)(1-\tau'(t)\rho){\rm d}\rho{\rm d}\sigma.
\end {align}
Now, let us suppose $\tau(t)\neq0$ and integrate by parts the last term in \eqref{5.28}. We get
\begin{align}\label{5.29}
&\int_{0}^{1}e^{-2\tau(t)\rho}u_{t}(x,t-\tau(t)\rho)u_{tt}(x,t-\tau(t)\rho)(1-\tau'(t)\rho){\rm d}\rho\nonumber\\
=&-\tau^{-3}(t)\int_{0}^{1}e^{-2\tau(t)\rho}u_{\rho}(x,t-\tau(t)\rho)u_{\rho\rho}(x,t-\tau(t)\rho)(1-\tau'(t)\rho){\rm d}\rho\nonumber\\
=&\tau^{-3}(t)\int_{0}^{1}e^{-2\tau(t)\rho}u_{\rho}(x,t-\tau(t)\rho)u_{\rho\rho}(x,t-\tau(t)\rho)(1-\tau'(t)\rho){\rm d}\rho\nonumber\\
&-\tau'(t)\tau^{-3}(t)\int_{0}^{1}e^{-2\tau(t)\rho}u_{\rho}^{2}(x,t-\tau(t)\rho){\rm d}\rho\nonumber\\
&-2\tau^{-2}(t)\int_{0}^{1}e^{-2\tau(t)\rho}u_{\rho}^{2}(x,t-\tau(t)\rho)(1-\tau'(t)\rho){\rm d}\rho\nonumber\\
&-\tau^{-3}(t)\left[e^{-2\tau(t)\rho}u_{\rho}^{2}(x,t-\tau(t)\rho)(1-\tau'(t)\rho)\right]_{0}^{1},
\end {align}
then we have
\begin{align}\label{5.30}
&\int_{0}^{1}e^{-2\tau(t)\rho}u_{t}(x,t-\tau(t)\rho)u_{tt}(x,t-\tau(t)\rho)(1-\tau'(t)\rho){\rm d}\rho\nonumber\\
=&-\frac{1}{2}\tau'(t)\tau^{-1}(t)\int_{0}^{1}e^{-2\tau(t)\rho}u_{t}^2(x,t-\tau(t)\rho){\rm d}\rho\nonumber\\
&-\int_{0}^{1}e^{-2\tau(t)\rho}u_{t}^2(x,t-\tau(t)\rho)(1-\tau'(t)\rho){\rm d}\rho\nonumber\\
&-\frac{\tau^{-1}(t)}{2}e^{-2\tau(t)}u_{t}^2(x,t-\tau(t))(1-\tau'(t))+\frac{\tau^{-1}(t)}{2}u_{t}^{2}(x,t).
\end {align}
Inserting \eqref{5.30} in \eqref{5.28}, we obtain
\begin{align}\label{5.31}
I'(t)=&\zeta\tau'(t)\int_{\Gamma_1}\int_{0}^{1}e^{-2\tau(t)\rho}u_{t}^{2}(x,t-\tau(t)\rho){\rm d}\rho{\rm d}\sigma\nonumber\\
&-2\zeta\tau'(t)\tau(t)\int_{\Gamma_1}\int_{0}^{1}e^{-2\tau(t)\rho}u_{t}^{2}(x,t-\tau(t)\rho){\rm d}\rho{\rm d}\sigma\nonumber\\
&-\zeta\tau'(t)\int_{\Gamma_1}\int_{0}^{1}e^{-2\tau(t)\rho}u_{t}^{2}(x,t-\tau(t)\rho){\rm d}\rho{\rm d}\sigma\nonumber\\
&-2\zeta\tau(t)\int_{\Gamma_1}\int_{0}^{1}e^{-2\tau(t)\rho}u_{t}^{2}(x,t-\tau(t)\rho)(1-\tau'(t)\rho){\rm d}\rho{\rm d}\sigma\nonumber\\
&-\zeta e^{-2\tau(t)}\int_{\Gamma_1}u_{t}^{2}(x,t-\tau(t)\rho)(1-\tau'(t)\rho){\rm d}\sigma+\zeta\int_{\Gamma_1}u_{t}^{2}(x,t){\rm d}\sigma\nonumber\\
=&-2\zeta\tau(t)\int_{\Gamma_1}\int_{0}^{1}e^{-2\tau(t)\rho}u_{t}^{2}(x,t-\tau(t)\rho){\rm d}\rho{\rm d}\sigma\nonumber\\
&-\zeta e^{-2\tau(t)}\int_{\Gamma_1}u_{t}^{2}(x,t-\tau(t)\rho)(1-\tau'(t)\rho){\rm d}\sigma+\zeta\int_{\Gamma_1}u_{t}^{2}(x,t){\rm d}\sigma,
\end {align}
from which immediately follows estimate \eqref{5.27} for $t$ such that $\tau(t)\neq0$. However, If  $\tau(t)=0$, from \eqref{5.14} we get
\begin{align}
I'(t)=&\zeta\tau'(t)\int_{\Gamma_1}\int_{0}^{1}e^{-2\tau(t)\rho}u_{t}^{2}(x,t-\tau(t)\rho){\rm d}\rho{\rm d}\sigma\nonumber\\
\leq&\zeta d \int_{\Gamma_1}\int_{0}^{1}u_{t}^{2}(x,t){\rm d}\rho{\rm d}\sigma\nonumber\\
=&\zeta d\int_{\Gamma_1}u_{t}^{2}(x,t){\rm d}\sigma\nonumber\\
=&\zeta d\int_{\Gamma_1}u_{t}^{2}(x,t){\rm d}\sigma-2I(t).\nonumber
\end {align}
Then, we obtain \eqref{5.27}.

{\bf Proof of Theorem \ref{th5.3}.}
From \eqref{5.3}, \eqref{5.14}, \eqref{5.18} and \eqref{5.27},  then from \eqref{5.8}, we get
\begin{align}\label{5.33}
L'(t)=&E'(t)+\varepsilon_1 G'(t)+\varepsilon_2 H'(t)+\varepsilon_3 I'(t)\nonumber\\
\leq&-\left[-\varepsilon_1-\varepsilon_2(\delta-g_0)\right]\|u_t(t)\|_2^2-\left[C-\varepsilon_1\left(1+\frac{\mu_{1}^{2}}{4\delta}\right)-\varepsilon_2\mu_1-\varepsilon_3\zeta\right]\|u_t(t)\|_{2,\Gamma_1}^2\nonumber\\
&-\left[\frac{1}{2}g(t)-\varepsilon_1\left(\delta-l+2c\delta \right)-\varepsilon_2\left(\delta+2\delta(1-l)^2\right)\right]\|\nabla u(t)\|_2^2\nonumber\\
&-\left(\alpha-\frac{\alpha\varepsilon_2}{2}\right)\|\nabla u_{t}(t)\|_{2}^{2}+\left[\frac{1}{2}-\frac{\varepsilon_2(1-l)}{4\delta\lambda_1}\right](g'\circ\nabla u)(t)\nonumber\\
&+\left[\frac{\varepsilon_1(1-l)}{4\delta\lambda_1}+\varepsilon_2\left(\frac{1-l}{4\delta}+\left(2\delta+\frac{1}{4\delta}\right)(1-l)^2+\frac{\alpha(1-l)}{2}
+\frac{1-l}{2\delta\lambda_1}\right)\right](g\circ\nabla u)(t)\nonumber\\
&-\left[C-\frac{\varepsilon_1\mu_{2}^2}{4\delta}-\varepsilon_2\mu_2\right]\|u_t(t-\tau(t))\|_{2,\Gamma_1}^{2}-2\varepsilon_3I(t),
\end {align}
By the trace inequality
\begin{align}
\|u_t\|_{2,\Gamma_1}^2\leq C\|u_t\|_{W^{1,2}(\Omega)} \leq C_1\|u_t\|_2^2,\nonumber
\end {align}
then we have
\begin{align}\label{5.34}
L'(t)=&E'(t)+\varepsilon_1 G'(t)+\varepsilon_2 H'(t)+\varepsilon_3 I'(t)\nonumber\\
\leq&-\left[-\varepsilon_1-\varepsilon_2(\delta-g_0)+C_1\left(\alpha-\frac{\alpha\varepsilon_2}{2}\right)\right]\|u_t(t)\|_2^2\nonumber\\
&-\left[\frac{1}{2}g(t)-\varepsilon_1\left(\delta-l+2c\delta \right)-\varepsilon_2\left(\delta+2\delta(1-l)^2\right)\right]\|\nabla u(t)\|_2^2\nonumber\\
&-\left[C-\varepsilon_1\left(1+\frac{\mu_{1}^{2}}{4\delta}\right)
-\varepsilon_2\mu_1-\varepsilon_3\xi\right]\|u_t(t)\|_{2,\Gamma_1}^2+\left[\frac{1}{2}-\frac{\varepsilon_2(1-l)}{4\delta\lambda_1}\right](g'\circ\nabla u)(t)\nonumber\\
&+\left[\frac{\varepsilon_1(1-l)}{4\delta\lambda_1}+\varepsilon_2\left(\frac{1-l}{4\delta}+\left(2\delta+\frac{1}{4\delta}\right)(1-l)^2+\frac{\alpha(1-l)}{2}
+\frac{1-l}{2\delta\lambda_1}\right)\right](g\circ\nabla u)(t)\nonumber\\
&-\left[C-\frac{\varepsilon_1\mu_{2}^2}{4\delta}-\varepsilon_2\mu_2\right]\|u_t(t-\tau(t))\|_{2,\Gamma_1}^{2}-2\varepsilon_3I(t).
\end {align}
At this point, we choose $\varepsilon_1$, $\varepsilon_2$ and $\varepsilon_3$ so small that \eqref{5.12}  remain valid and
\begin{align}
k_1=-\varepsilon_1-\varepsilon_2(\delta-g_0)+C_1\left(\alpha-\frac{\alpha\varepsilon_2}{2}\right)>0,\nonumber
\end {align}
\begin{align}
k_2=\frac{1}{2}g(t)-\varepsilon_1\left(\delta-l+2c\delta \right)-\varepsilon_2\left(\delta+2\delta(1-l)^2\right)>0,\nonumber
\end {align}
\begin{align}
k_3=C-\varepsilon_1\left(1+\frac{\mu_{1}^{2}}{4\delta}\right)
-\varepsilon_2\mu_1-\varepsilon_3\xi>0,\nonumber
\end {align}
\begin{align}
k_4=\frac{1}{2}-\frac{\varepsilon_2(1-l)}{4\delta\lambda_1}>0,\nonumber
\end {align}
\begin{align}
k_5=\frac{\varepsilon_1(1-l)}{4\delta\lambda_1}+\varepsilon_2\left(\frac{1-l}{4\delta}+\left(2\delta+\frac{1}{4\delta}\right)(1-l)^2+\frac{\alpha(1-l)}{2}
+\frac{1-l}{2\delta\lambda_1}\right),\nonumber
\end {align}
\begin{align}
k_6=C-\frac{\varepsilon_1\mu_{2}^2}{4\delta}-\varepsilon_2\mu_2>0.\nonumber
\end {align}
Therefore, \eqref{5.34} takes the form
\begin{align}\label{5.35}
L'(t)\leq&-k_1\|u_t(t)\|_2^2-k_2\|\nabla u(t)\|_2^2-k_3\|u_{t}(t)\|_{2,\Gamma_1}^{2}+k_4(g'\circ\nabla u)(t)+k_5\left(g\circ\nabla u\right)(t)\nonumber\\
&-k_6\|u_t(t-\tau(t))\|_{2,\Gamma_1}^{2}-2\varepsilon_3I(t).
\end {align}
Since $I(t)\geq0$ and by \eqref{2.2}, $(G2)$ there exists a positive constant $M$ such that
\begin{align}\label{5.36}
L'(t)\leq-ME(t)+k_5(g\circ\nabla u)(t),\ \forall t\geq t_0.
\end {align}
Multiplying \eqref{5.36} by $\xi(t)$, we have
\begin{align}\label{5.37}
\xi(t)L'(t)\leq-M\xi(t)E(t)+k_5\xi(t)(g\circ\nabla u)(t),\ \forall t\geq t_0,
\end {align}
Because $\xi$ and $g$ are nonincreasing, we get
\begin{align}
\xi(t)\int_{0}^{t}g(t-s)\|\nabla u(t)-\nabla u(s)\|_{2}^{2}ds\leq&-\int_{0}^{t}g'(t-s)\|\nabla u(t)-\nabla u(s)\|_{2}^{2}ds\nonumber\\
\leq&-2E'(t)\nonumber
\end {align}
Inserting the last inequality in \eqref{5.37}, we obtain
\begin{align}\label{5.38}
\xi(t)L'(t)+2k_5E'(t)\leq-M\xi(t)E(t),\ \forall t\geq t_0.
\end {align}
Now,we define
\begin{align}
H(t)=\xi L(t) +2k_5E(t).\nonumber
\end {align}
Since  $\xi(t)$ is nonincreasing positive function, we can easily get that $H\sim E$. Thus \eqref{5.38} implies that
\begin{align}
H'(t)\leq -k\xi(t)H(t), \ \forall t\geq t_0 ,\nonumber
\end {align}
for some $k>0$. Then, by direct integration over $(t_0, t)$, we have
\begin{align}
H(t)\leq H(t_0)e^{ -k\int_{t_0}^t\xi(s){\rm d}s}, \ \forall t\geq t_0. \nonumber
\end {align}
Consequently, using the equivalent relations of $H(t)$ and $E(t)$ , we can conclude
\begin{align}
E(t)\leq k_5\Phi(t_0)e^{ -k\int_{t_0}^t\xi(s){\rm d}s}=Ke^{-k\int_{t_0}^t\xi(s){\rm d}s}, \ \forall t\geq t_0, \nonumber
\end {align}
where $k_5$ is a positive constant and $K=k_5\Phi(t_0)$.
This completes the proof.

\end{document}